\documentclass{ifacconf}

\usepackage{filecontents}

\makeatletter
\let\old@ssect\@ssect 
\makeatother

\usepackage{natbib}
\usepackage{hyperref}

\makeatletter
\def\@ssect#1#2#3#4#5#6{%
  \NR@gettitle{#6}
  \old@ssect{#1}{#2}{#3}{#4}{#5}{#6}
}
\makeatother

\usepackage{graphicx}      
\graphicspath{{Figures/}}  
\usepackage{float}         
\usepackage{placeins}

\usepackage{amsmath}

\usepackage{natbib}        

\usepackage{booktabs}
\usepackage{makecell}
\usepackage{tabularx}
\usepackage{multirow}

\usepackage{hyperref}

\begin{document}
\begin{frontmatter}

\title{Impact of Spatial and Technology Aggregation on Optimal Energy System Design}

\author[a]{Shruthi Patil}  
\author[a]{Leander Kotzur} 
\author[b]{Detlef Stolten}

\address[a]{Institute of Techno-economic Systems Analysis (IEK-3), Forschungszentrum Jülich GmbH, Wilhelm-Johnen-Str., D 52428, Germany} 
\address[b]{Chair for Fuel Cells, RWTH Aachen University, c/o Institute of Techno-economic Systems Analysis (IEK-3), Forschungszentrum Jülich GmbH, Wilhelm-Johnen-Str., D 52428, Germany}

\begin{abstract}
Designing an optimal energy system with large shares of renewable energy sources is computationally challenging. Considering greater spatial horizon and level of detail, during the design, exacerbates this challenge. 

This paper investigates spatial and technology aggregation of energy system model, as a complexity-reduction technique. To that end, a novel two-step aggregation scheme based on model parameters such as Variable Renewable Energy Sources (VRES) time series and capacities, transmission capacities and distances, etc, is introduced. First, model regions are spatially aggregated to obtain a reduced region set. The aggregation is based on a holistic approach that considers all the model parameters and spatial contiguity of the regions. Next, technology aggregation is performed on each VRES, present in each newly-defined region. Each VRES is aggregated based on the temporal profiles to obtain a representative set. The impact of these aggregations on the accuracy and computational complexity of a cost-optimal energy system design is analyzed for a European energy system scenario. 

The aggregations are performed to obtain different combinations of number of regions and VRES types, and the results are benchmarked against an initial spatial resolution of 96 regions and 68 VRES types in each region. The results show that the system costs deviate significantly when lower number of regions and/or VRES types are considered. As the spatial resolution is increased in terms of both number of regions and VRES types, the system cost fluctuates at first and stabilizes at some point, approaching the benchmark value. Optimal combination can be determined based on an acceptable cost deviation and the point of stabilization. For instance, if $<5\%$ deviation is acceptable, 33 regions and 38 VRES types in each region is optimal. With this setting, the system cost is under-estimated by $4.42\%$ but the run time is reduced by $92.95\%$.

\end{abstract}


\end{frontmatter}


\section{Introduction} \label{sec:Introduction}

\subsection{Background: Spatio-temporal Energy System \\Optimization Models}

In the light of international agreements to reduce Greenhouse Gas (GHG) emissions (\cite{agreement2015paris}), participating governments have introduced policies that promote increasing the shares of Renewable Energy Sources (RES) in the energy systems. However, the design and real-time operation of an energy system, with large shares of RES, is very challenging because of the intermittent nature of primary energy sources (e.g., wind and solar energy). The time duration during which electricity demand is high does not match the time during which electricity can be harnessed from RES. In addition, energy demand sites might sometimes be far away from the locations eligible for installation of different RES (\cite{samsatli2015general}). In short, there exists a spatio-temporal gap between electricity demand and supply. To reduce this gap, transmission options to transport the excess
output of one location to nearby locations and storage options to store the surplus energy for later use are necessary. In addition, conversion technologies to convert produced electricity into a storable form like hydrogen, and vice versa are required. The combination of these technologies make up a very complex energy system.

The deployment and operation of such a complex system is not straightforward. A long-term planning and strategic deployment of various technologies that are employed in generation, conversion, storage, and transport of electricity is required. In order to do so, it is vital to assess the impact of different decisions relating to the size, location, combination, operation rate, etc. of these technologies. 

An Energy System Optimization Model (ESOM) that accounts for the spatial and temporal dependencies of these technologies and the dynamic nature of the energy demand, can be employed to support the decision-making process. These spatio-temporal ESOMs minimize a certain objective by optimizing different technologies' location, capacity, and their operation, subject to various system-specific and user-defined constraints (\cite{decarolis2017formalizing}).

In general, ESOMs are formulated in the following manner:

Given an ESOM, which contains:
\begin{itemize}
	\item Spatial description: Geographical area and the network of
	regions within it. Each region in this network is treated
    as a single node with connections to neighboring regions for transmission of surplus electricity, under the "copper plate"
	assumption (\cite{cao2018incorporating}). 

	\item A set of technologies within each region: 
	\begin{itemize}
		\item Different generation, storage, conversion, and transport 
		technologies
		\item Minimum and maximum capacity of each technology
		\item Capital, operating, and maintenance costs, etc.
	\end{itemize}
\end{itemize}

Subject to constraints such as:
\begin{itemize}
	\item The energy demand 
	\item Resource availability and its maximum operation limit
\end{itemize}

Determine:
\begin{itemize}
	\item Size and location of different technologies
	\item Operation rate of each technology
\end{itemize}

\begin{figure*}
    \begin{center}
      \includegraphics[width=\textwidth, scale=1]{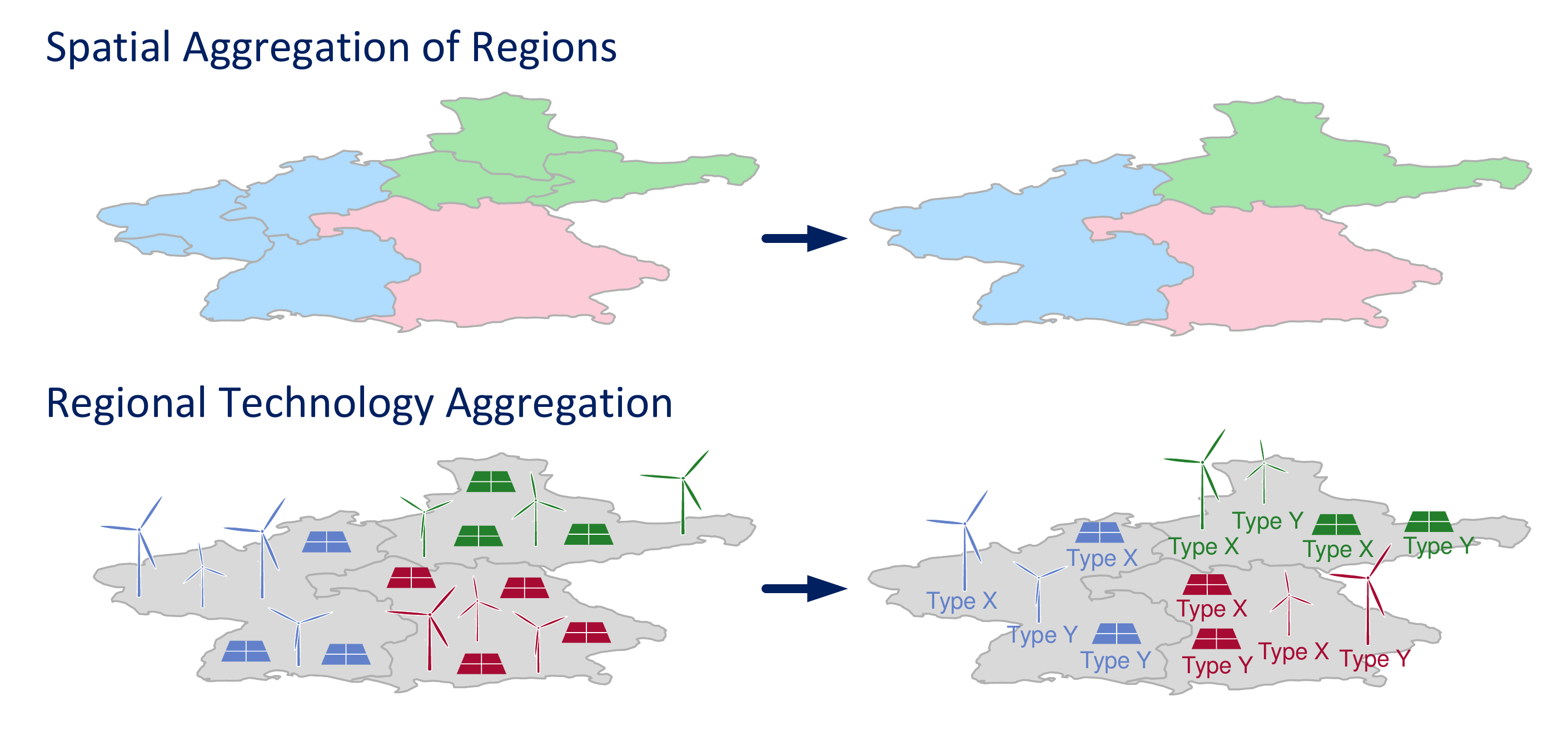}
  
      \caption{Pictorial description of spatial aggregation of regions and technology aggregation within each newly-defined region.}
        \label{fig:spagat_basic_depiction}
    \end{center}
\end{figure*}

Such that a certain performance criterion is optimized, such as, minimization of the total cost of the energy system or minimization of GHG emissions.

Although different spatio-temporal ESOMs that are developed have the same general formulation as described above, they differ mainly in the optimization criteria. For instance, \cite{welder2018spatio} developed an ESOM whose optimization criteria is to reduce the Total Annual Cost (TAC) of the energy system. In their work, \cite{samsatli2018multi} developed a multi-objective ESOM. The ESOM can be optimized to obtain an energy system with minimal cost, maximal profit, minimal CO2 emissions, or maximal energy production, or a desired combination of these. 

\subsection{Data Aggregation for Complexity Reduction}

One of the major challenges concerning ESOMs is their associated computational complexity (\cite{pfenninger2014energy}). According to \cite{ridha2020complexity}, ESOMs have four complexity dimensions. They are:
\begin{enumerate}
    \item Mathematical complexity: The mathematical formulation of the model i.e., linear, mixed-integer linear, or non-linear, etc. Also, the ability to take stochastic behavior of the system into consideration. 
    \item Temporal complexity: The temporal resolution and horizon of the model.
    \item Spatial complexity: The spatial resolution and horizon of the model.
    \item System scope: The system's parts and level of detail that are taken into account in the model.
\end{enumerate}

A steady improvement in the quality and availability of data allows for the incorporation of greater details in ESOMs (\cite{priesmann2019complex}). However, as the level of detail increases so does the complexity along one, a combination, or all of the above-mentioned complexity dimensions. Increase in the complexity of ESOMs necessitates more computational power and higher solving times. Beyond a certain level of complexity, ESOMs run the risk of becoming computationally intractable (\cite{frew2016temporal}).

A popular approach to reduce the complexity of ESOM is to reduce its size by employing data aggregation techniques as a pre-processing step, prior to optimization. These data aggregation techniques could be broadly classified into three types, namely, temporal, spatial, and technology aggregation. 

\subsubsection{Temporal aggregation:}

The basic idea behind temporal aggregation is to coarsen the temporal resolution of highly resolved demand and supply time series, thereby reducing the temporal complexity of the ESOM (\cite{cao2019classification}). Fundamentally, there are two ways to perform temporal aggregation - (i) directly reducing the number of time steps and (ii) reducing the number of periods. Direct reduction of time steps can be achieved either by down-sampling (i.e., the time series are divided into slices of fixed size and the values within each slice are averaged) or by segmentation (i.e., adjacent time steps are merged if there exits mutual similarity between them). On the other hand, reducing the number of periods is based on the understanding that there might exist similarities between different periods in the time series and not just between adjacent time steps. Here, the time series are sliced into periods (e.g. each day) and periods with similar profiles are grouped to form typical periods (\cite{kotzur2020modeler}).

Some works apply these techniques in their original form (\cite{pandvzic2014effect, nahmmacher2016carpe, kotzur2018impact}), while more recent works either (i) improve upon these techniques. For e.g. \cite{de2020variable} maintain the resolution during critical periods and down-sample the rest of the time series or (ii) use a combination of these techniques. For e.g., \cite{fazlollahi2014multi} first group the time series into typical days (period length is a day) and within each typical day, segmentation is performed. A review of various temporal aggregation techniques can be found in \cite{hoffmann2020review}.

\subsubsection{Spatial aggregation:}
As the name suggests, spatial aggregation reduces the ESOM data along the spatial dimension. Here, model regions are spatially aggregated by merging contiguous regions with similar properties (\cite{grubesic2014spatial}). This reduces the spatial resolution of the ESOM, thereby reducing its spatial complexity. 

\subsubsection{Technology aggregation:}
Technology aggregation involves aggregation of technologies based on similar properties. E.g., all wind turbines present in a region can be aggregated based on the similarity in their temporal profiles to obtain a representative set. In terms of the above-mentioned complexity dimensions, it reduces the system scope by reducing the level of detail considered. A pictorial description of spatial and technology aggregation is shown in Figure \ref{fig:spagat_basic_depiction}.

\subsection{Objective and Structure of the Paper}

This paper focuses on spatial and technology aggregation of ESOM data, as a complexity-reduction technique. We begin with a highly resolved ESOM instance, and:
\begin{enumerate}
    \item First, the model regions are spatially aggregated. 
    \item Next, the Variable Renewable Energy Sources (VRES) (i.e., wind turbine and photovoltaic) profiles are aggregated, within each newly-defined region.
\end{enumerate}

The general objective here is to reduce the size of the ESOM instance without compromising on the accuracy of its optimization results. In other words, the aggregation of ESOM should strike a balance between its accuracy and its computational complexity. In order to achieve this, either the number of aggregated regions can be kept considerably low, while keeping the number of VRES types per aggregated region reasonably high, or vice versa. In the light of this, the main objective of this paper is to determine how many regions and how many VRES types within each region would be sufficient to strike a balance between accuracy and computational complexity of ESOMs.

The remainder of this paper is structured as follows: Section \ref{sec:StateofResearch} explores the previous work on spatial and technology aggregation and highlights the research gaps. In Section \ref{sec:Methodology}, the methodology employed to close the highlighted research gaps is introduced. The experimental results are found in Section \ref{sec:Results}. A summary and discussion can be found in Section \ref{sec:SummaryandDiscussion}. The main conclusions of the study are given in Section \ref{sec:Conclusion}. Finally, some future areas of research are discussed in Section \ref{sec:Outlook}.     

\section{State of Research} \label{sec:StateofResearch}

\subsection{Spatial Aggregation of Regions}

According to \cite{fischer1980regional}, aggregation of regions involves forming of homogeneous region sets from an initial set of regions. Each homogeneous region set consists of spatially contiguous regions (i.e., neighboring regions) which show high degree of similarity with respect to an attribute or a set of attributes. In essence, aggregation of regions is a spatially constrained clustering problem. Methods developed to solve this problem include - (i) sequentially applying a conventional clustering technique with no regard to geography and then grouping identified similar regions only if they are contiguous, (ii) adding the x and y coordinates of each region's centroid as its two additional attributes, and (iii) explicitly including the spatial contiguity constraint in the clustering procedure (\cite{duque2012max}). Alternatively, graph theory-based algorithms can be used. They reduce a connected graph into connected sub-graphs, maximizing a similarity criterion within each sub-graph.   

In the context of energy systems analysis, spatial aggregation is popularly viewed as a network reduction problem. For e.g., in their work, \cite{horsch2017role} consider the European electricity transmission network and apply k-means clustering (\cite{lloyd1982least}) to reduce this network. The aim here is to reduce the network while maintaining the major transmission corridors. The attributes considered for similarity definition are load and generation capacities. In addition to these attributes, geographical coordinates are considered to ensure spatial contiguity. The authors note three ways in which their work could be improved - (i) considering and comparing different clustering algorithms, (ii) considering higher number of clusters than the 362 clusters considered, and (iii) considering additional attributes.

According to \cite{biener2020grid}, it is important to consider electrical grid characteristics during network reduction as it ensures accurate grid representation in the reduced ESOM. To that end, they introduce a network reduction method that takes electrical distances between nodes as the similarity defining attribute. A combination of density-based graph clustering (\cite{zhou2009graph}) and agglomerative hierarchical clustering (\cite{ward1963hierarchical}) methods are used to cluster the nodes. Using quality indicators such as root-mean-square error, they compare their method with the one developed by \cite{horsch2017role}. Most of the quality indicators chosen show that their method outperforms the one developed by \cite{horsch2017role}.

\cite{cao2018incorporating} start with the German transmission grid. They consider the marginal costs of the total power supply as the similarity defining attribute and spectral clustering (\cite{fiedler1973algebraic}), a graph theory-based algorithm, is used to cluster the nodes. The number of clusters is set to 20 and the ESOM is run for both the reduced case and the fully resolved case. In comparison to the fully resolved case, the reduced case shows a deviation of $7.4\%$  in the optimization results, but the computing time is drastically reduced to $4.3\%$. The authors note that finding the optimal number of clusters could be a future research topic.

The HotMaps Horizon2020 project (\cite{scaramuzzino2019integrated}) aims to promote the development of transnational 
renewable energy policies and strategies. In order to show that regions belonging to different nations are similar,
the European NUTS3 regions (\cite{eurostat1995nomenclature}) are aggregated, based not only on various energy potential indicators (e.g., wind, solar, agricultural residues potentials, etc.) but also economic (e.g., electricity and gas prices, etc.) and socio-demographic (e.g., population, GDP, etc.) attributes. In addition to these attributes, geographic locations are also considered. k-means clustering is used to cluster the regions and the NbClust tool (\cite{malika2014nbclust}) is used to identify the optimal number of clusters. The results show that NUTS3 regions can be reduced to 17 clusters 
using this method. It is noteworthy that although the geographic locations of the regions are considered during clustering, 
not all regions in the clusters are contiguous. 

As a part of the e-Highway2050 project, \cite{anderski2015european} aggregate NUTS3 regions. Various attributes are considered to define similarity between regions- population, mean wind speed, mean solar irradiation, thermal installed capacity, hydro installed capacity, and agricultural areas and natural grasslands. These attributes are weighted depending on their significance. In addition, geographic locations are considered. The clustering method employed is from the Python module ClusterPy (\cite{ClusterPy}), which is based on k-means and tabu-search algorithms. The algorithm is run on one country at a time in order to avoid grouping regions of 
different countries. The decision regarding the number of groups in each country is based on its total area, population, and load. The method yields 105 clusters with contiguous regions within each cluster. After consulting the transmission system operators, these clustered regions are further reduced to 96 regions.

As opposed to the above works, \cite{siala2019impact} begin with a high resolution raster data. They introduce a novel spatial aggregation algorithm which is based on k-means++ (\cite{vassilvitskii2006k}) and max-p regions (\cite{duque2012max}) algorithms. Aggregation of the data cells is based on wind potential, photovoltaic potential, or electricity demand at a time. An ESOM is run for each of these cases and the results are compared with each other and also for the case of national borders. They conclude that 
region definition based on any one of the above-mentioned characteristics leads to better optimization results compared to national borders. As a future work, they suggest combining various characteristics during spatial aggregation.

\subsection{Technology Aggregation}

In the literature, it is commonplace to simply aggregate each ESOM component data, within each defined region (\cite{kotzur2020modeler, cao2019classification, siala2019impact}). For e.g., demand and generation time series are averaged and their capacities are summed (\cite{welder2018spatio}). However, some works investigate the clustering of time series applied to energy systems analysis. 

Clustering of demand time series is seen in some publications (\cite{de2017comparing}, \cite{rasanen2009feature}, \cite{sun2016c}). 
With respect to generation time series, \cite{joubert2016optimisation} optimize wind farm locations using mean-variance portfolio optimization method. In a pre-processing step, the wind farms are clustered using agglomerative hierarchical clustering. The optimization is run for different number of clusters and the results are compared with those obtained in the case of optimization run with unclustered data, in order to determine the optimal number of clusters. The authors conclude that even in the case of optimal number of clusters, there is a marginal deviation in the optimization results compared to the unclustered solution. However, they note that clustering the data has the benefit of reduced computation time. 

With an aim to find the most suitable clustering technique to cluster photovoltaic power time series, \cite{munshi2016photovoltaic} apply various clustering approaches ranging from conventional techniques such as k-means and hierarchical clustering to genetic algorithms (\cite{goldberg2006genetic}) like ant colony and bat clustering. The authors conclude that bat clustering exhibits the best performance, but is computationally intensive. 

In the context of ESOMs, \cite{caglayan2021robust} address the aggregation of VRES. Here, the time series of each VRES are grouped based on their respective Levelized Cost of Electricity (LCOE) and the time series within each group are averaged. This procedure is repeated for all defined regions. It is seen that considering more than one time series per VRES, per region leads to lower TAC compared to one time series per VRES, per region. Increasing the VRES resolution leads to the availability of more cost-competitive locations to choose from, thereby bringing the overall system costs down. 

\cite{radu2021model} introduce a two-stage procedure to identify and keep only the most relevant VRES locations and discard those that have little impact on the results of optimization. In the first stage, a simplified version of the ESOM is run with full set of VRES locations and the locations chosen during optimization are deemed relevant. In the second stage, the full version of ESOM is run with these VRES locations. The performance is evaluated by comparing the results with those obtained when full version of ESOM is run with full set of VRES locations. The results show that more than $90\%$ of the relevant VRES locations are correctly identified by the procedure and the memory consumption and solver time are reduced by up to $41\%$ and $46\%$, respectively. 

\cite{frysztacki2021strong} extend the work of \cite{horsch2017role} to investigate the individual and combined effects of network resolution and VRES resolution on the results of ESOM. To that end, they consider three cases - (i) the transmission network is clustered and the VRES sites and other system components are aggregated to the nearest network node (ii) the county-level network resolution is maintained and the VRES resolution in each country is varied, and (iii) both network and VRES resolutions are varied. The results show that system costs are under-estimated at low network resolutions, as network bottlenecks are not revealed at lower resolutions. On the other hand, considering low VRES resolution over-estimates system costs due to unavailability of cost-competitive locations. The authors conclude that both network and VRES resolution should be sufficiently high in order to accurately estimate the optimal system costs.

\subsection{Research Gaps}

Previous works on spatial aggregation focus on either one or a few of its parameters. A holistic approach, considering all the model parameters is missing. If homogeneous region groups are formed based on a similarity definition that takes into account all the parameters, aggregation of each parameter, within each group based on simple techniques such as mean or sum would be justified. 

In order to ensure spatial contiguity in the resulting region groups, it is commonplace to consider geographic locations as additional attributes during grouping. This, however, does not always ensure spatial contiguity as seen in \cite{scaramuzzino2019integrated}. This is especially true if several attributes are considered during grouping and all attributes are weighted equally. Therefore, explicitly including the spatial contiguity constraints during region grouping could be more suitable. 

Finally, owing to the high variance in the time series of VRES such as wind turbines and photovoltaics, it would be relevant to consider a high spatial resolution of VRES, per newly-defined region. Since this would lead to an increase in computational complexity, we propose to cluster each VRES, in each region based on the similarity in temporal profiles, to obtain a representative set. 

\begin{figure*}
    \begin{center}
      \includegraphics[width=\textwidth]{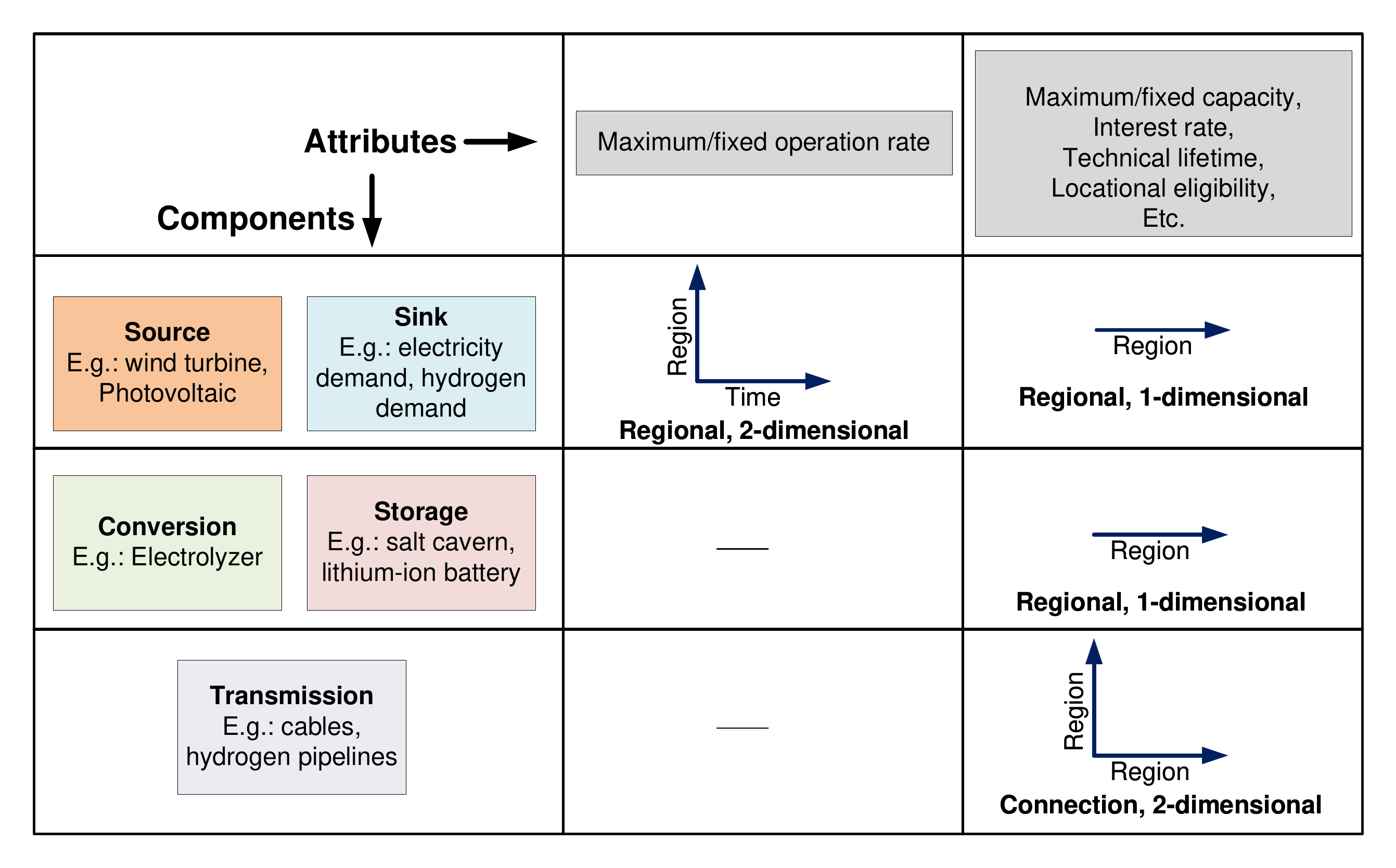}
  
      \caption{The general data structure of an Energy System Optimization Model (ESOM) instance.}
        \label{fig:components_attributes_datastructure}
    \end{center}
\end{figure*}

\section{Methodology} \label{sec:Methodology}

This section is divided into five subsections. The first subsection provides details about the ESOM employed and its setup. The second subsection provides the general approach followed. The next two subsections introduce the approaches adopted for spatial and technology aggregation, respectively. The last subsection describes the experimental design. 

\subsection{Energy System Optimization Model Details}

For the analysis, an open source optimization framework called Framework for Integrated Energy System Assessment (FINE) (\cite{welder2018spatio}) is employed. Parameterising FINE essentially involves adding various energy system components with their corresponding data. Figure \ref{fig:components_attributes_datastructure} shows various component classes, their components, and the data attributes corresponding to each component. These attributes can be classified into the following two types:

\begin{enumerate}
    \item \textit{Regional attributes}: These attributes are region specific. Further, 
    they can be 1-dimensional (region) or 2-dimensional (region * time). Examples of 1-dimensional and 2-dimensional regional attributes are maximum capacity and the capacity factor time series of wind turbine, respectively. 
    \item \textit{Connection attributes}: These attributes characterize the connections between region pairs. They are always 2-dimensional (region * region). An example of a connection attribute is capacity of a DC cable between two regions. 
\end{enumerate}

Several attributes belonging to various system components, with varying dimensions and data types make up a very complex data structure. This data is stored as a netCDF file. Python's xarray module (\cite{hoyer2017xarray}) is used to read in the saved netCDF files. The computations are performed on the read in xarray dataset. 

\subsection{General Approach}

\begin{figure*}
    \begin{center}
      \includegraphics[width=\textwidth]{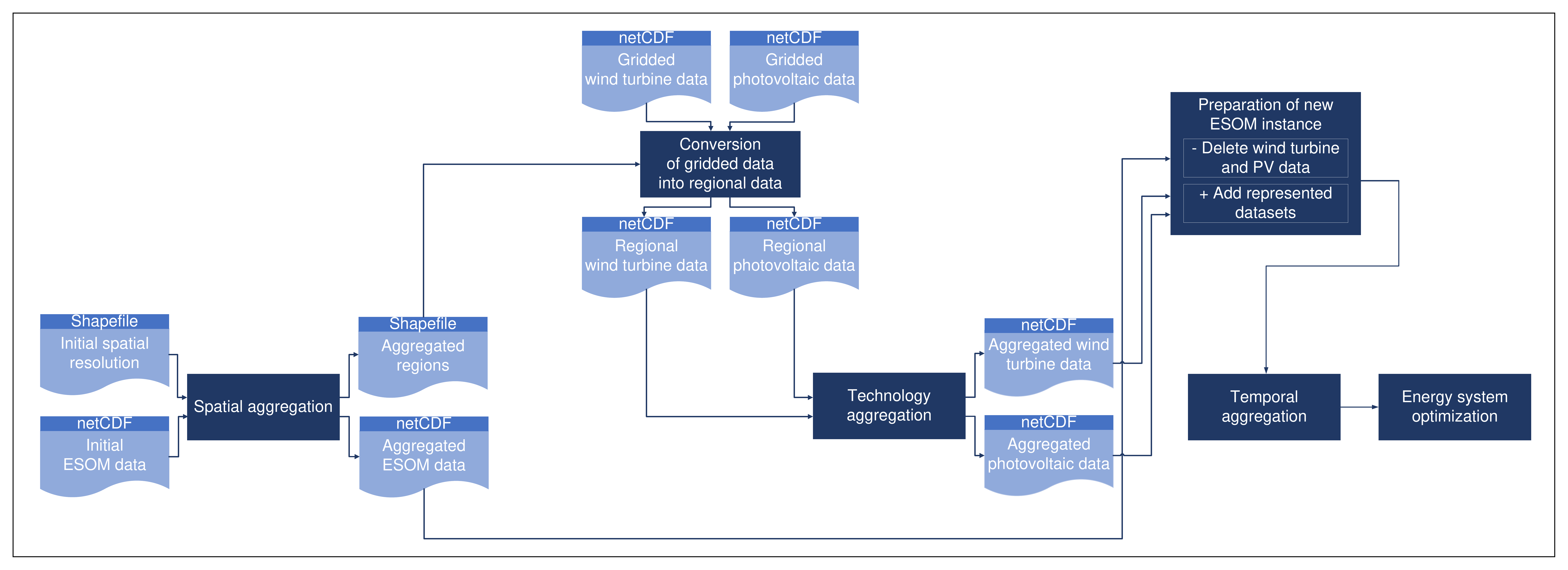}
  
      \caption{Block diagram depicting the general workflow adopted for the experiments.}
        \label{fig:block_diagram_general_approach}
    \end{center}
\end{figure*}

Figure \ref{fig:block_diagram_general_approach} shows the general workflow adopted for the experiments. 
Initially, a FINE instance is set up by adding various system components. In this FINE instance, just one type of each component is present in every region. For e.g., there is one photovoltaic in each region whose capacity factor (CF) time series and corresponding capacity are an aggregation of all simulated photovoltaics for that region. Spatial aggregation takes as its input this
FINE instance (in the form of a netCDF file) and a shapefile with the initial spatial resolution. It results in two outputs - a netCDF file and a shapefile with aggregated data and spatial resolution, respectively. 

The open source tools Renewable Energy Simulation toolkit (RESKit) (\cite{ryberg2019reskit}) and Geospatial Land Availability for Energy Systems (GLAES) (\cite{Ryberg2018}) are employed to place and simulate the wind turbines and photovoltaics. A gridded data with 50 * 50 spatial resolution across Europe is obtained using these tools. The resulting datasets can be overlapped with a shapefile to convert the gridded data into regional data. This regional data is the input for technology aggregation block. Aggregation is performed separately for wind turbines and photovoltaics.  

The wind turbines and photovoltaics present in the result of grouping are now deleted and the represented sets are added. A new FINE instance is set up based on this aggregated data. 

In order to keep the computation time within reasonable limits, temporal aggregation is performed on the resulting dataset. Finally, optimization is performed on the spatially, technologically, and temporally aggregated FINE instance. 

\subsection{Spatial Aggregation}

\subsubsection{Algorithm:} In order to aggregate the regions the Hess model (\cite{hess1965nonpartisan}), also known as k-medoids clustering, is employed here. The aim is to partition a given region set V to form k groups. 

If the number of regions in V is $n$ and $i$ and $j$ are two arbitrary regions, then Hess model uses the 
following ${n}^2$ binary variables:

\begin{align*}
  x_{ij}=\begin{cases}
    1, & \text{if $i$ is assigned to a group with center at $j$}\\
    0, & \text{otherwise}.
  \end{cases}
\end{align*}

The Hess model is formulated in the following manner:

\begin{subequations}

\begin{equation}\label{eq:obj}
    min \sum_{i \in V} \sum_{j \in V} D(i,j) x_{ij}
\end{equation}

where, $D(i,j)$ is the distance between the regions $i$ and $j$.

subject to the constraints:
\begin{equation}\label{eq:c1}
     \sum_{j \in V} x_{ij} = 1  \qquad	\forall i \in V
\end{equation}

\begin{equation}\label{eq:c2}
     \sum_{j \in V} x_{jj} = k
\end{equation}

\begin{equation}\label{eq:c3}
    x_{ij} \leq x_{jj}  \qquad	\forall i,j \in V
\end{equation}

\end{subequations}

Constraints \ref{eq:c1} ensure that each region is assigned to a group and only one group. Constraint \ref{eq:c2} ensures that k groups are formed. Finally, constraint \ref{eq:c3} ensures that a region $i$ is assigned to region $j$ only if $j$ is chosen to be a group's center. 

In order to ensure spatial contiguity in the resulting region groups the contiguity constraints, first introduced by \cite{oehrlein2017cutting} for spatial aggregation, are employed here. The constraint formulation is based on the concept of $(a,b)-separators$. An $(a,b)-separator$ can be defined as a subset of regions $C \subseteq V \backslash \{a,b\}$ that, if removed, would destroy all paths connecting regions $a$ and $b$. Since removing all regions present in C would lead to disconnected regions $a$ and $b$, it is obvious that for regions $a$ and $b$ to belong to a group, at least one of the regions from C should also be present in the group (\cite{validi2020imposing}). Mathematically, this condition is expressed as follows:

\begin{equation}\label{eq:c4}
    \sum_{c \in C} x_{cb} \geq x_{ab}  \qquad	\forall C \subseteq V \backslash \{a,b\}
\end{equation}

With an aim to speed up the calculations, the model is first solved without the contiguity constraints. The resulting region groups are scouted for disconnected region pairs and subsequently, contiguity constraint for these region pairs are added and the model is solved again. This iterative process is stopped once all the regions in each group are connected. \footnote{The method introduced here is published in the Python package  \href{https://github.com/FZJ-IEK3-VSA/tsam/}{tsam - Time Series Aggregation Module} and can be easily applied}

\subsubsection{Distance measure:} Considering the complexity of the dataset, a custom distance is defined to determine the distance between region pairs. This distance definition is based on residual sum of squares and works on the values normalized across each data attribute. Mathematically, the distance between two regions \textit{a} and \textit{b} is defined as:

\begin{equation}\label{eq:customDistance}
    D(a, b) = D_{r\_1d}(a, b) + D_{r\_2d}(a, b) + D_{c\_2d}(a, b)
\end{equation}

where,

$D_{r\_1d}(a, b)$ is the cumulative distance of all 1-dimensional regional attributes between the two regions, and is defined by: 

\begin{equation}\label{eq:dist_reg1d}
    D_{r\_1d}(a, b) = \sum_{i \in D_{r\_1d}} {(i(a) - i(b))}^2
\end{equation}

$D_{r\_2d}(a, b)$ is the cumulative distance of all 2-dimensional regional attributes between the two regions, and is defined by: 

\begin{equation}\label{eq:dist_reg2d}
    D_{r\_2d}(a, b) = \sum_{i \in D_{r\_2d}} \sum_{t=1}^{t=T} {(i(a, t) - i(b, t))}^2
\end{equation}

Here, t = 1, 2,...,T are the time steps

And, $D_{c\_2d}(a, b)$ is the cumulative distance of all 2-dimensional connection attributes between the two regions, and is defined by: 

\begin{equation}\label{eq:dist_conn2d}
    D_{c\_2d}(a, b) = \sum_{i \in D_{c\_2d}} {(1 - i(a,b))}^2
\end{equation}

Connection attributes indicate how strongly two regions are connected and their normalized values lie in the range [0,1]. Therefore, these values are converted into a distance meaning by subtracting from 1 as shown in Equation \ref{eq:dist_conn2d}.

\subsubsection{Connectivity matrix:} 
A connectivity matrix indicates which region pairs are connected and is employed in the algorithm described above to ensure spatial contiguity. In this matrix, the value corresponding to a region pair is 1 if they are spatially contiguous, otherwise 0. Two regions are deemed contiguous if:

\begin{enumerate}
    \item the regions' borders touch at least at one point
    \item one of the regions is an island and the other its nearest mainland region 
    \item there is a transmission line or a pipeline running between the two regions
\end{enumerate} 

\subsubsection{Data aggregation:}

Once the new region set is obtained, the data within each region group is aggregated. The aggregation method varies depending on the attribute. The aggregation method corresponding to each attribute is shown in Table \ref{tab:attributerepresentation}. 

\begin{table}[h]
	\caption{Aggregation method employed for each data attribute}
	\begin{center}
	\begin{tabularx}{\columnwidth}{X X}
		\toprule
		\multicolumn{1}{c}{\thead{\textbf{Attribute}}} & \multicolumn{1}{c}{\thead{\textbf{Aggregation method}}} \\
		\toprule
		
		Maximum operation rate & Weighted mean (weights being its corresponding maximum capacity)\\ 
		Fixed operation rate & Sum\\ 
		Maximum capacity & Sum \\ 
		Fixed capacity & Sum \\ 
		Locational eligibility & Boolean OR \\ 
        Investment per capacity & Mean \\ 
        Investment if built & Boolean OR \\ 
        Opex per operation & Mean \\ 
        Opex per capacity & Mean \\ 
        Opex if built &  Boolean OR \\ 
        Interest rate &  Mean \\ 
        Economic lifetime &  Mean \\
        Losses  & Mean \\ 
        Distances &  Mean \\ 
        Commodity cost &  Mean \\ 
        Commodity revenue &  Mean \\ 
        Opex per charge operation &  Mean \\ 
        Opex per discharge operation  & Mean  \\ 
        Technical lifetime &  Sum \\ 
        Reactances & Sum \\
        
        \bottomrule
	\end{tabularx}
    \end{center}
	\label{tab:attributerepresentation}
\end{table}

\subsection{Technology Aggregation}

To aggregate the VRES in each region, Python Scikit-learn's agglomerative hierarchical clustering is used. Here, the connectivity matrix is not provided, as the exact locations of the technologies, within a region, is irrelevant under the copper plate assumption. 

\begin{figure}[H]%
\begin{center}
\includegraphics[width=8.4cm]{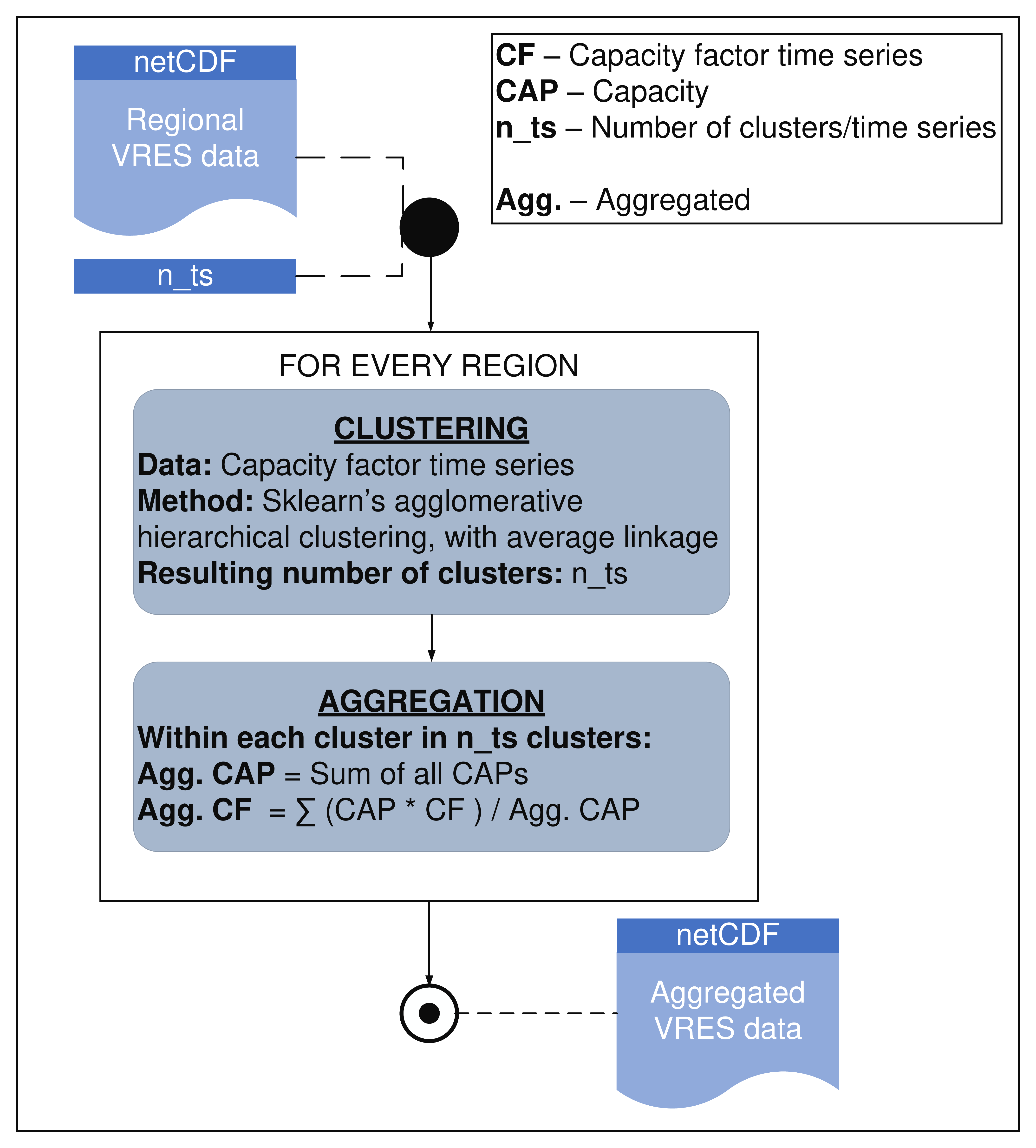}    
\caption{Flowchart of technology aggregation algorithm} 
\label{fig:representation_workflow}
\end{center}
\end{figure}

Figure \ref{fig:representation_workflow} describes the implemented technology representation algorithm. The regional data (Capacity Factor (CF) time series and capacities) of a particular VRES and the desired number ($n_{ts}$) of time series per region are input. The algorithm works on one region at a time. First, the clustering technique is run. It clusters the CF time series to obtain $n_{ts}$ clusters. Next, aggregation of data is performed. Within each cluster, the aggregated capacity is determined by the sum of all capacities belonging to that cluster. The aggregated CF time series is obtained by taking the weighted mean of all the CF time series belonging to the cluster, with the respective capacities being their weights. 

\subsection{Experimental Design}

\subsubsection{Spatial and temporal scope and resolution of ESOM:}

\begin{enumerate}
    \item \textbf{Spatial scope of the ESOM:} A European energy system scenario (\cite{caglayan2021robust}) is considered in this paper. A FINE instance is set up for the same. 
    
    \item \textbf{Initial spatial resolution of ESOM:} The region definition suggested in the e-highway study is considered. In this study, the geographical area of Europe is divided into 96 regions.
    
    \item \textbf{Temporal scope of ESOM:} The data of 1 year is considered. 

    \item \textbf{Temporal resolution of ESOM:} The dataset has hourly temporal resolution, with 8760 time steps for 1 year. Prior to optimization, temporal aggregation is performed. The resolution is reduced to 40 typical days with 8 segments within each typical day. For this purpose, the method developed by \cite{hoffmann2020review} is employed. 
    
    In order to accurately access the impact of spatial and technology aggregation on the optimization results, the effects of temporal aggregation should be nullified across all experimental runs. Therefore, the temporal aggregation is run once for the highest spatial resolution and the number of VRES time series. The resulting cluster order is saved and in all the successive runs, the data is temporally aggregated  to obtain the same cluster order. 
    
\end{enumerate}

\subsubsection{Evaluation Method:} 

Initially, the combination of highest spatial resolution (i.e., initial spatial resolution of 96 regions) and number of VRES time series present in the smallest region (i.e., 68) are chosen and the general procedure shown in Figure \ref{fig:block_diagram_general_approach} is followed. The results of this run forms the benchmark. Successively, various combinations of number of aggregated regions and VRES time series per region are chosen and the procedure is repeated. 

The complexity indicator, in each case, is the total time taken and the accuracy indicator is the optimization objective i.e., TAC of all the components considered in the FINE instance. 

\section{Results} \label{sec:Results}

\subsection{Spatial Aggregation}

\begin{figure*}[!htbp]
    \begin{center}
      \includegraphics[width=\textwidth]{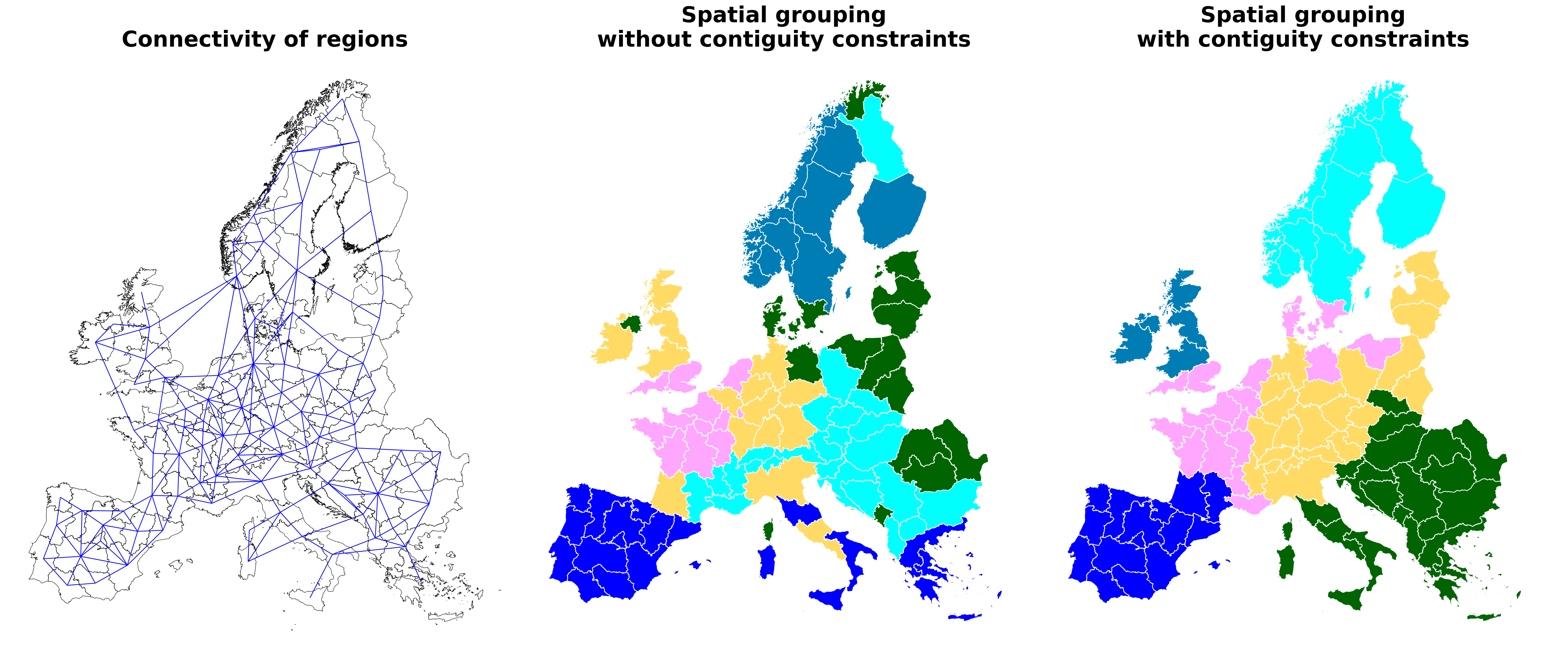}
  
      \caption{[Left] The initial 96 e-highway regions. The blue lines show connectivity between regions, as per the obtained connectivity matrix. [Middle] The 6 region groups obtained when spatial aggregation is performed without contiguity constraints. [Right] The 6 region groups obtained when spatial aggregation is performed with contiguity constraints.}
        \label{fig:grouping6reg_with_and_without_contiguity}
    \end{center}
\end{figure*}

For the purpose of comparison, the initial 96 regions are aggregated to obtain 6 region groups without and with the contiguity constrains. The results are shown in Figure \ref{fig:grouping6reg_with_and_without_contiguity} along with the connectivity of regions, as per the obtained connectivity matrix. 

When aggregation is performed without contiguity constraints, the resulting region groups are not fully connected. For e.g., the north-most region in Norway is only connected to its immediate neighbors in Norway and Finland. However, it is grouped with some southern regions. A similar observation can be made in other region groups too. 

\begin{figure*}[!htbp]
    \begin{center}
      \includegraphics[width=\textwidth]{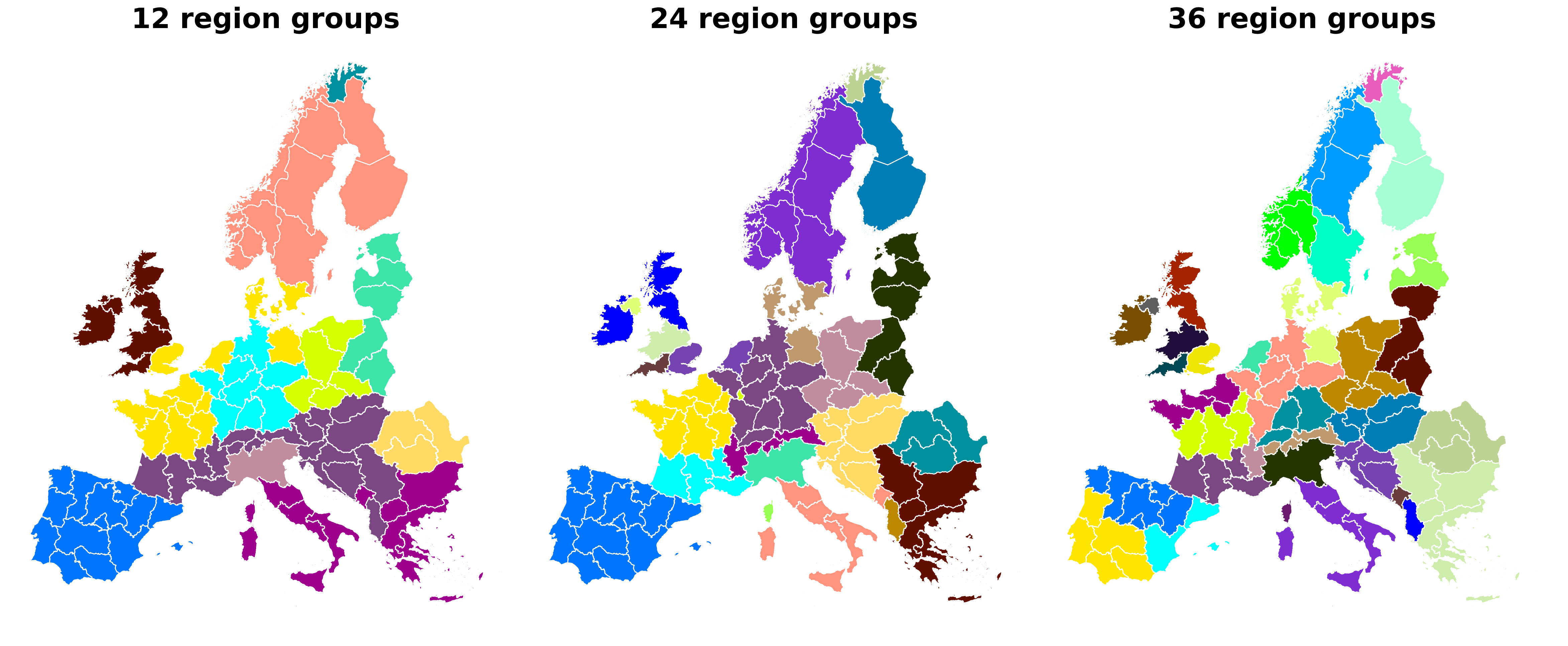}
  
      \caption{Region groups obtained when the 96 regions are spatially aggregated with contiguity constraints, to obtain 12, 24, and 36 region groups.}
        \label{fig:grouping_12_24_36reg}
    \end{center}
\end{figure*}

The plot on the right in Figure \ref{fig:grouping6reg_with_and_without_contiguity} shows that aggregating the regions with contiguity constraints ensures formation of region groups that are fully connected. Here, all the regions are compact except the one shown in pink. This region has some fragmented parts. Nonetheless, the original regions present in this group are connected, hence contiguity constraints are not violated here. Further, in Figure \ref{fig:grouping_12_24_36reg} it is seen that as the number of region groups is increased, the region groups get more compact.

\begin{figure*}[!htbp]
    \begin{center}
      \includegraphics[width=\textwidth, scale=0.5]{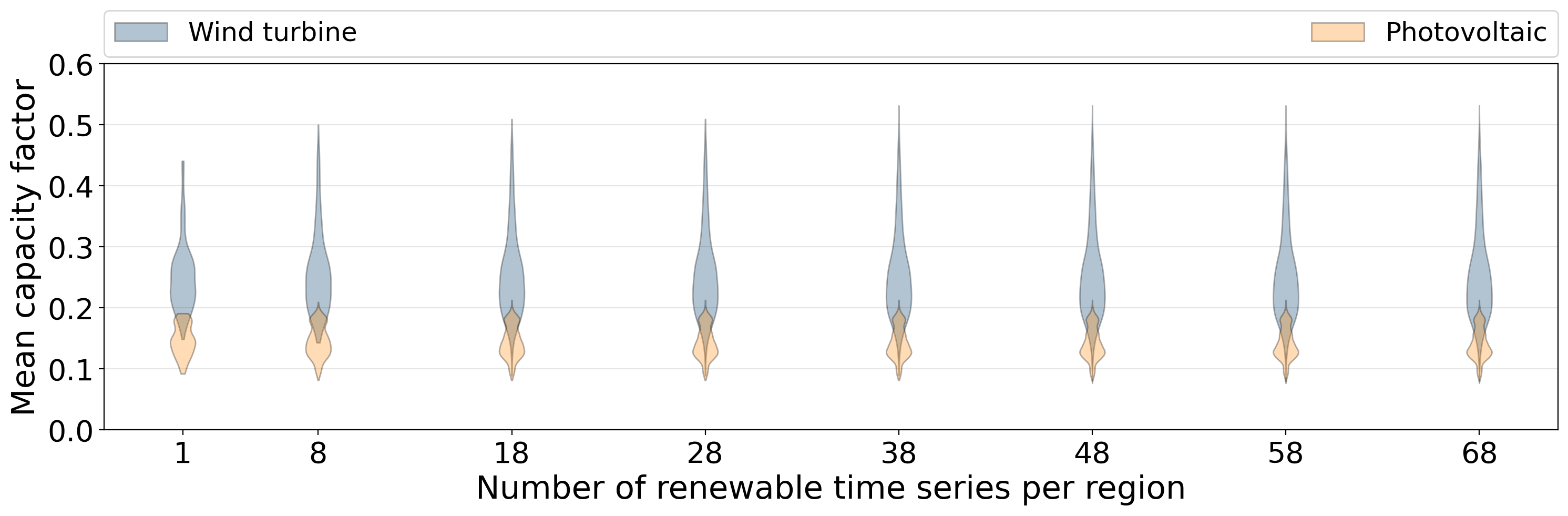}
  
      \caption{Distribution of installable capacities in all 96 e-highway regions, for different number of VRES time series per region.}
        \label{fig:V_installable}
    \end{center}
\end{figure*}

\subsection{Technology Aggregation}

For the 96 regions, the technology aggregation algorithm is run to obtain different number of VRES time series per region. The resulting capacity distribution is seen in Figure \ref{fig:V_installable}. Along the x-axis, the curves show the capacity distribution for different number of time series per region. The y-axis indicates the mean CF of each time series. The width of the curve, corresponding to each y-value, shows the capacity corresponding to the time series.

As the number of time series is increased, two changes can be noticed in wind turbine and photovoltaic distributions. The wind turbine data distribution narrows at mid-range i.e., at mean CF ranging between 0.2 and 0.3 and the photovoltaic data distribution narrows at mean CF ranging between 0.1 and 0.2. In both the distributions, at lower number of time series the extreme values are not captured well. As the number of time series is increased, the extreme values appear. At 38 time series, all the extreme values seem to be captured, as an increase in the number of time series further does not change the distribution of both photovoltaic and wind turbine data. 

\subsection{Impact of Spatial and Technology Aggregation on Optimization Results}

\begin{figure*}[!htbp]
    \begin{center}
      \includegraphics[width=\textwidth, scale=0.5]{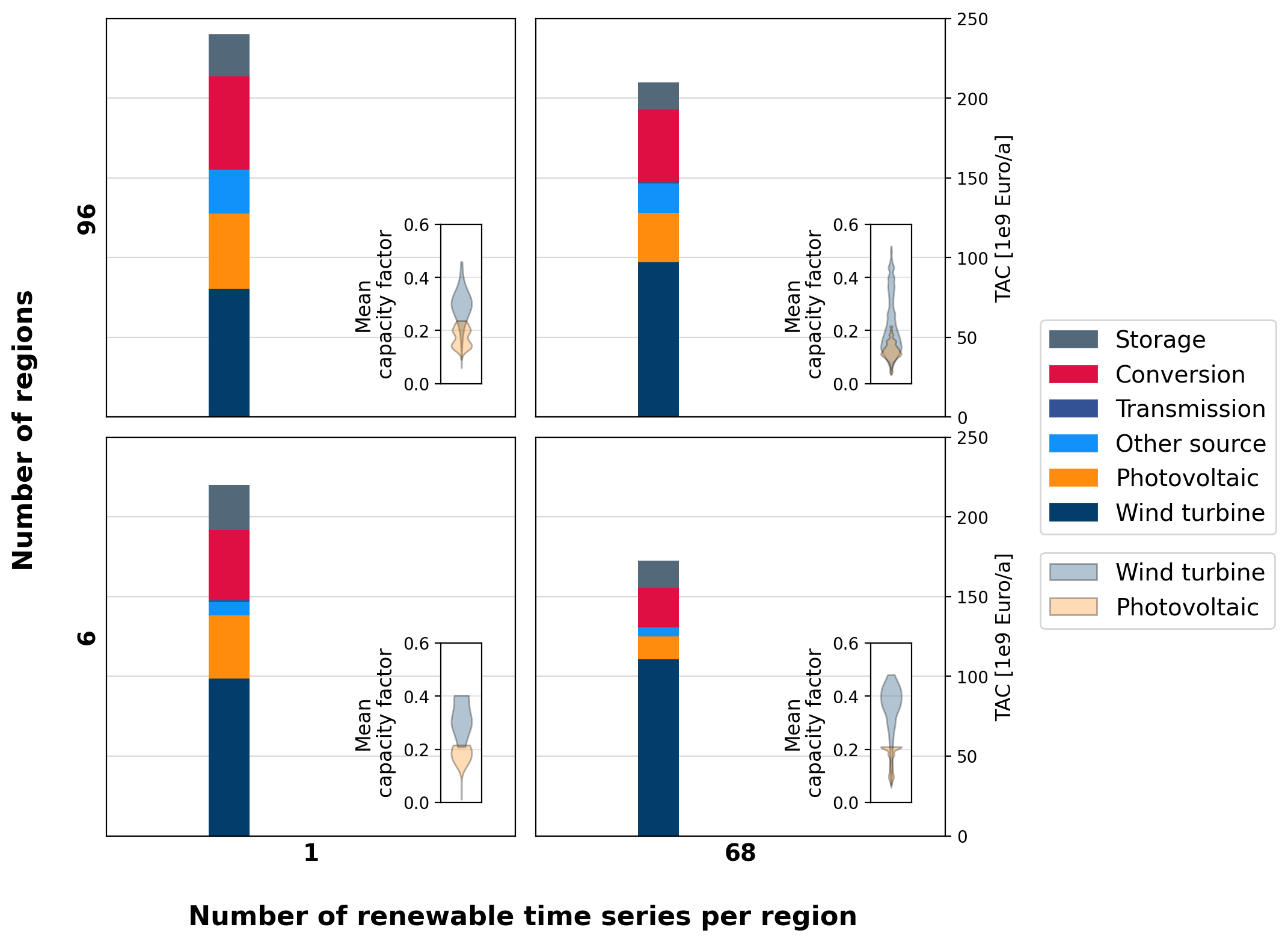}
  
      \caption{The results of optimization for the extreme parameter combinations (i.e., 6 and 96 regions and 1 and 68 VRES time series per region). For each of these combinations, TAC for different technologies is seen on the left and distribution of installed VRES capacities is seen on the right.}
        \label{fig:TAC_matrix_plus_violins}
    \end{center}
\end{figure*}

The lowest and the highest number of regions considered in this study are 6 and 96, respectively. The lowest and the highest number of VRES time series per region considered are 1 and 68, respectively. Initially, the optimization results obtained for these extreme combinations are analysed. Figure \ref{fig:TAC_matrix_plus_violins} shows the optimization results for these parameter combinations. For each parameter setting, a bar plot on the left shows the TAC for different technologies and on the right, the distribution of installed VRES capacities (i.e., distribution of optimal capacity and operational time series chosen during optimization) is seen. 

The top-right cell in this figure shows the benchmark results i.e., the results for 96 regions with 68 VRES time series per region. With this setting, the overall TAC is approximately 210 billion Euro/annum. The wind turbine's capacity distribution curve starts with a needle-like shape around 0.52 mean CF and slowly widens for lower mean CFs and narrows again at mean CF around 0.1. The photovoltaic's capacity distribution curve also has a similar shape. It starts at around 0.2 mean CF and slowly widens and narrows again at mean CF around 0.1.

Now, the top-left cell shows the optimization results when the number of regions is kept the same but the number of VRES time series is reduced to 1 per region. It can be observed that reducing the number of time series, while keeping the number of regions constant, leads to an increase in the TAC. The capacity distributions of VRES help explain this behaviour. It is seen in Figure \ref{fig:V_installable} that the extreme values are not captured well when 1 time series per region is considered. Due to the fact that wind turbines with high mean CFs are not available to choose from, most wind turbines that are installed have mean CFs around 0.3. In comparison to the benchmark, fewer wind turbines are installed bringing the TAC of wind turbines down. As an alternative, more photovoltaics and other source technologies are installed, thereby increasing their TACs. Increase in the installation of these technologies requires an increase in the installation of conversion and storage technologies. Therefore, the TAC of these technologies also increases.

In the bottom-right cell, the optimization results obtained when the number of regions is reduced to 6 but the number of VRES time series per region is kept the same as benchmark can be seen. It is observed that reducing the number of regions, while keeping the number of VRES time series per region constant, leads to a decrease in the TAC. Since each region is considered as a single node under the copper plate assumption, the size of these regions does not play a role during optimization. In other words, it does not matter if these regions are spread across Europe or just a country. The spatial details within these regions is also limited since each component is aggregated within these regions. These factors give it an effect that the energy system is small with a network of 6 regions, leading to an under-estimation of the overall TAC. 

The wind turbine's capacity distribution shows that the peaks are not captured well in this setting thereby installing more wind turbines with mean CFs around 0.4. When good locations for wind turbine installation is found in each region, owing to the size of these regions and the copper plate assumption, there is lesser need for alternative sources. Therefore, in comparison to the benchmark, more wind turbines and fewer photovoltaics and other sources are installed here. Owing to the decrease in the installation of these technologies, a decrease in the installation of conversion and storage technologies is observed. 

Finally, the bottom-left cell shows the effect of reducing both the number of regions and the number of VRES time series per region. Both, over-estimation of the TAC due to the decrease in the number of time series and the under-estimation due to the decrease in the number of regions is seen here. The overall TAC is approximately 220 billion Euro/annum. It comes close to the benchmark value. However, individual TACs and associated capacity distributions differ. 

The wind turbine's capacity distribution shows that only those with mean CFs between 0.2 and 0.4 are installed in this setting. Due to the unavailability of better locations for wind turbines, more photovoltaics with mean CFs around 0.2 are installed. Other sources are installed less when compared to the benchmark, leading to the decrease in the installation of conversion technologies. On the other hand, increase in photovoltaic installation has lead to the increase in installation of storage technologies. 

\begin{figure*}[!htbp]
    \begin{center}
      \includegraphics[width=\textwidth, scale=0.5]{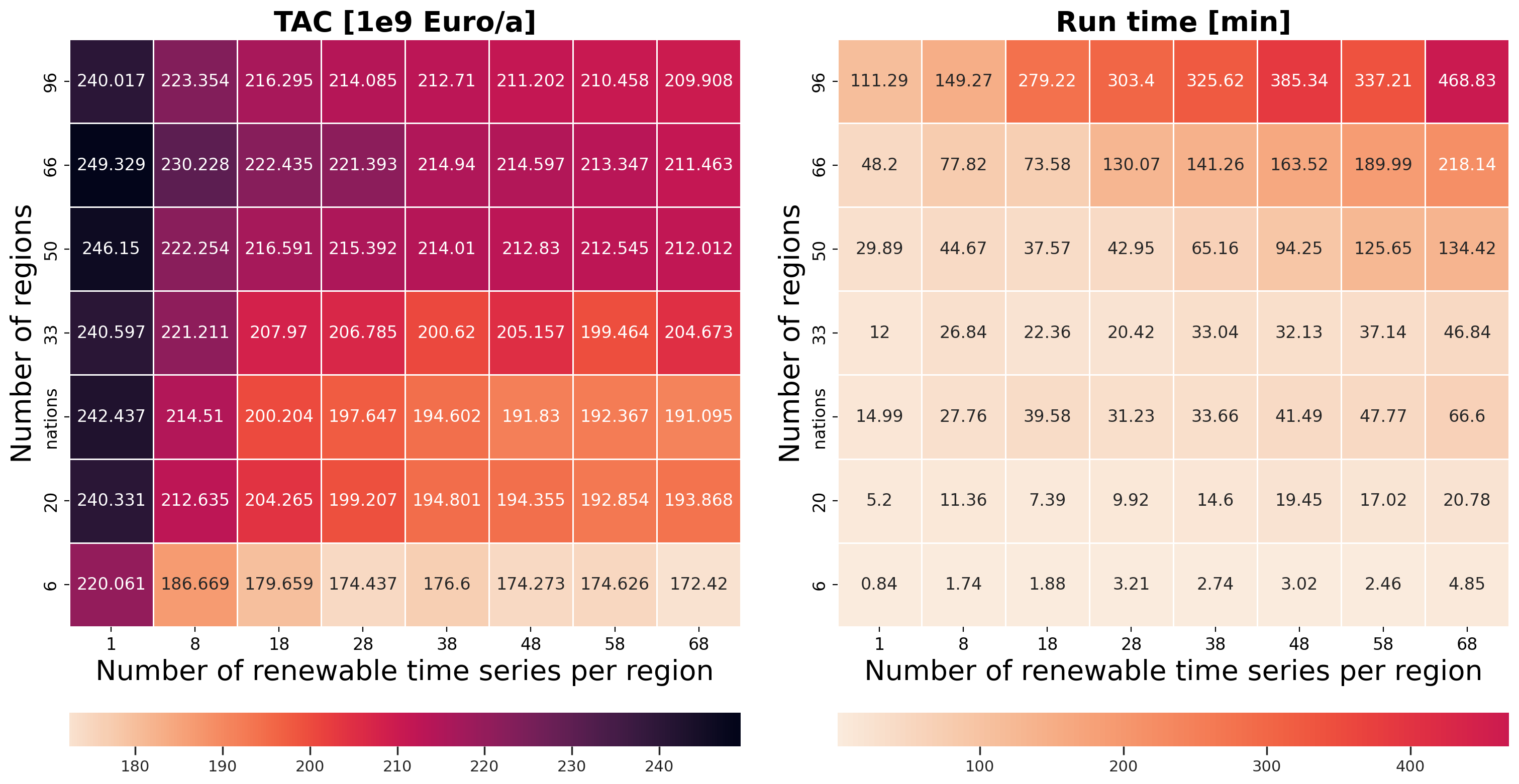}
  
      \caption{TAC obtained and the associated run time for each combination of parameters considered in this study.}
        \label{fig:TAC_time}
    \end{center}
\end{figure*}

Figure \ref{fig:TAC_time} shows the TAC obtained and the associated run time for each parameter setting. The general observation made earlier - decreasing the number of regions leads to an decrease in the TAC and reducing the number of VRES time series per region, increases the TAC - holds true overall. The run time is high for the benchmark setting and reduces in both the directions. 

It is noteworthy that in both the TAC and run time matrices, some exceptions to the general observations are visible. For instance, in the case of 33 aggregated regions reducing the time series from 58 to 48  leads to an increase in the TAC, but reducing further to 38 time series decreases the TAC. Also, the run time for 96 regions and 58 time series is 337.21 minutes. However, when the model is optimized for 96 regions and 48 time series the run time increases to 385.24 minutes. The reason for such exceptions could be that it is more challenging to find a global minimum in case of certain region groups and the set of representative VRES time series within these groups.

In addition to the different number of aggregated regions, the results for 33 national regions can also be seen in Figure \ref{fig:TAC_time}. Comparing the TACs obtained for 33 nations and 33 aggregated regions with the benchmark values, it can be observed that the TACs obtained for 33 nations deviate further from the benchmark. This shows that considering administrative boarders while designing an optimal energy system does not yield the best results. 

\begin{figure*}[!htbp]
    \begin{center}
      \includegraphics[width=\textwidth, scale=0.5]{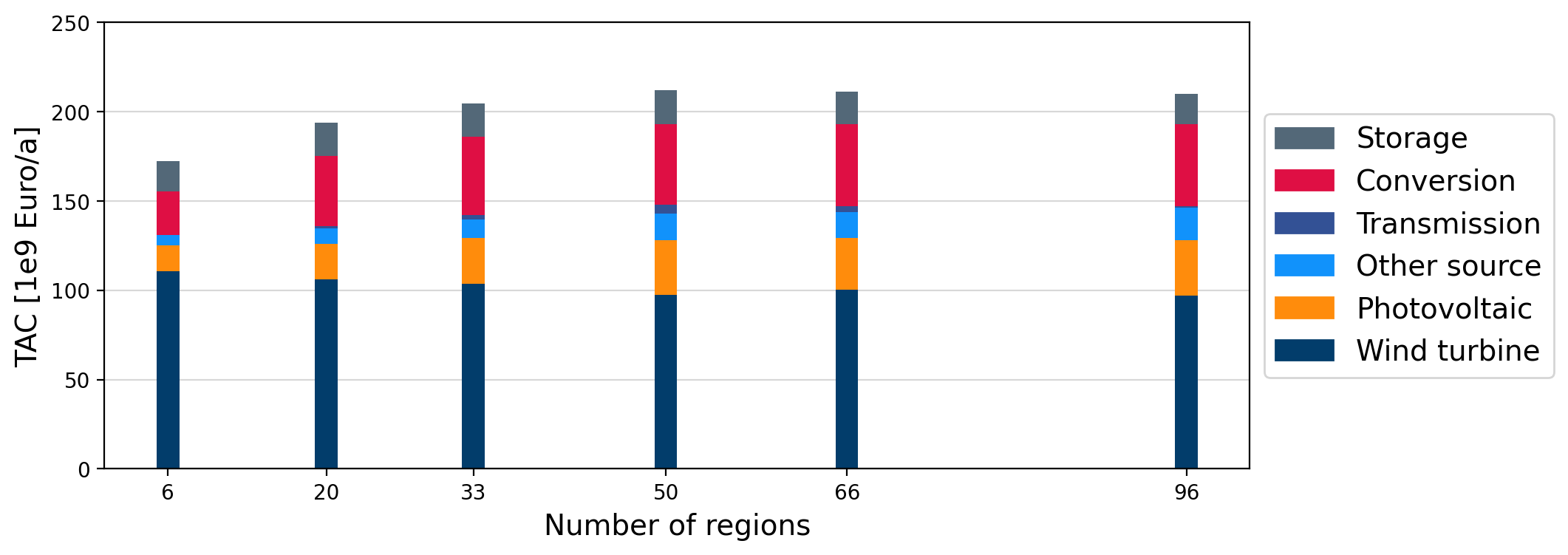}
  
      \caption{TAC obtained for different components in the case of differed number of aggregated regions, while considering 68 VRES time series per region.}
        \label{fig:TAC_for_diff_reg}
    \end{center}
\end{figure*}

\begin{figure*}[!htbp]
    \begin{center}
      \includegraphics[width=\textwidth, scale=0.5]{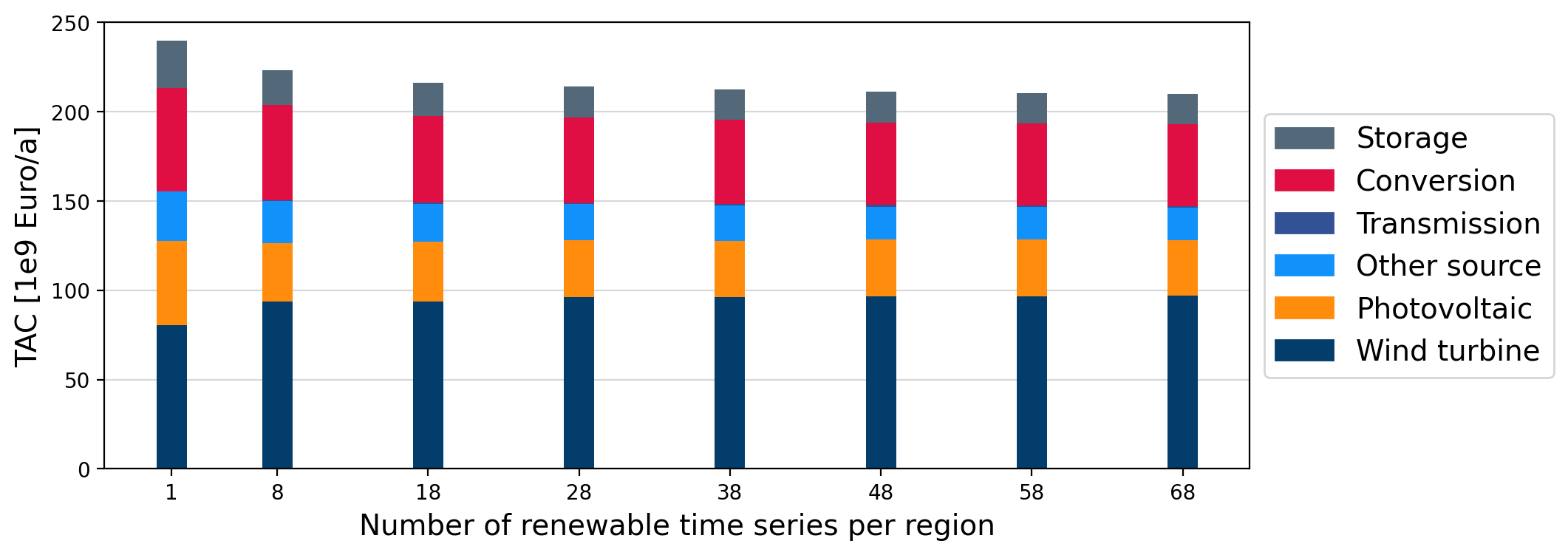}
  
      \caption{TAC obtained for different components in the case of 96 regions, when different number of aggregated VRES times series per region are considered}
        \label{fig:TAC_for_diff_ts}
    \end{center}
\end{figure*}

Figure \ref{fig:TAC_for_diff_reg} shows the TAC for different technologies when the number of regions is varied while considering 68 VRES time series within each region, in each case. It can be observed that as the number of regions is increased, the overall TAC and also the TAC of individual technologies stabilize. The optimal number of regions here can be determined using the elbow criterion i.e., the number of regions where a further increase would not change the TACs significantly. Similarly, in Figure \ref{fig:TAC_for_diff_ts} where the number of time series is varied for 96 regions, the optimal number of time series per region would be the number of time series where the TACs stabilize. 

\begin{figure*}[!htbp]
    \begin{center}
      \includegraphics[width=\textwidth, scale=0.5]{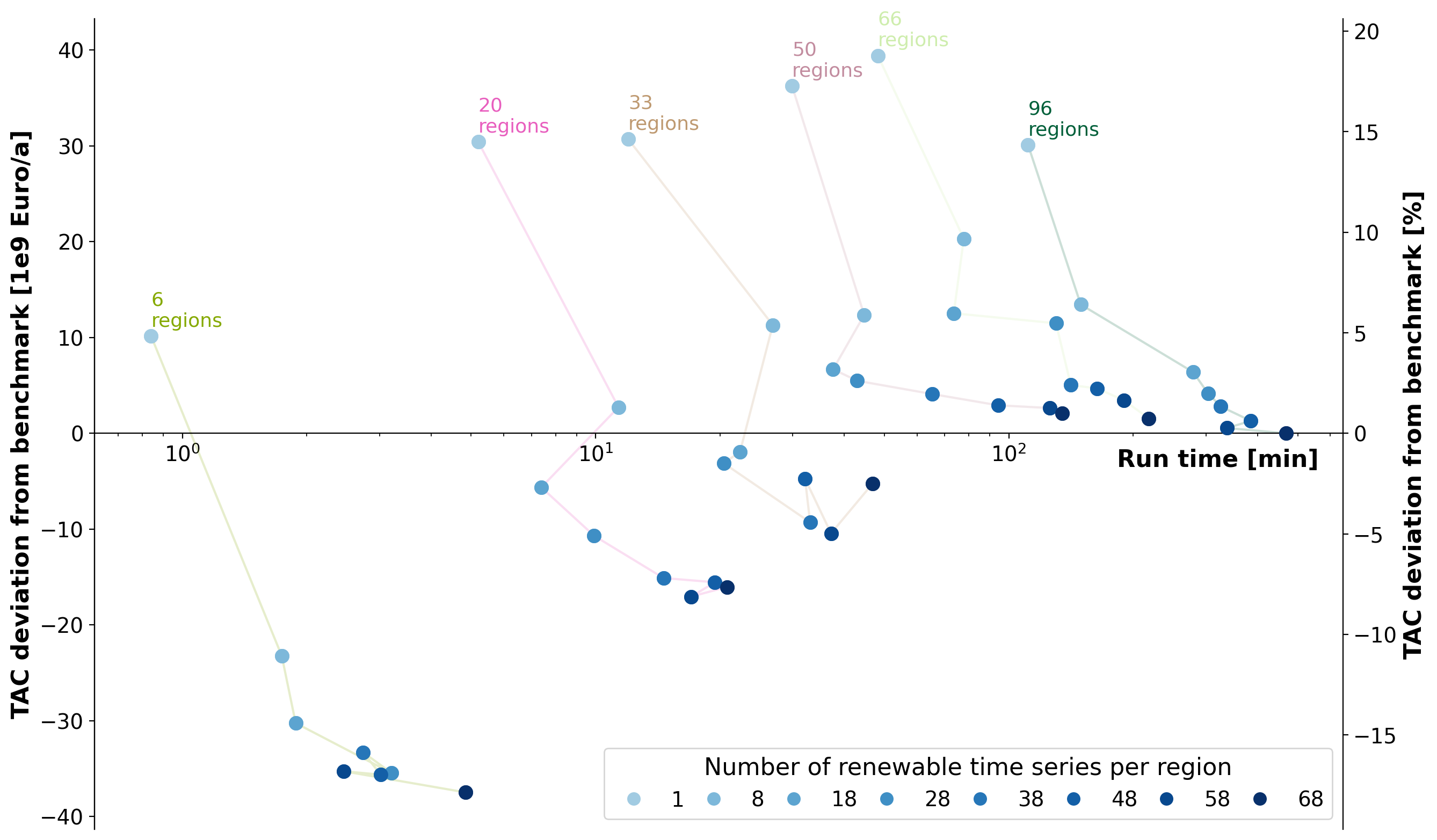}
  
      \caption{TAC deviation from the benchmark (represented on y-axis in billion Euro/annum on the left and in $\%$ on the right) versus the run time (represented on a logarithmic scale on x-axis) for each parameter setting. Each dot represents a particular parameter setting. The darker the shade of blue, the more the number of VRES time series considered and the lines connecting the dots represent a particular number of regions.}
        \label{fig:TAC_deviation_vs_time}
    \end{center}
\end{figure*}

Now, in order to determine the optimal combination of the number of aggregated regions and VRES time series within each region, it is necessary to apply the elbow criterion in both the directions of increase. In Figure \ref{fig:TAC_deviation_vs_time}, the run time (represented on a logarithmic scale on x-axis) and the TAC deviation from the benchmark (represented on y-axis in billion Euro/annum on the left and in $\%$ on the right), for each parameter setting can be seen. Each dot here represents a particular parameter setting. The darker the shade of blue, the more the number of VRES time series considered and the lines connecting the dots represent a particular number of regions. 

From the figure, it can be observed that for each regions set, as the number of VRES time series is increased the deviation is reducing drastically at first and seems to stabilise after a certain number of time series. On the other hand, keeping the number of VRES time series per region constant, if the number of regions is increased the deviation approaches 0, in each case.

\begin{figure}[H]%
\begin{center}
\includegraphics[width=8.4cm]{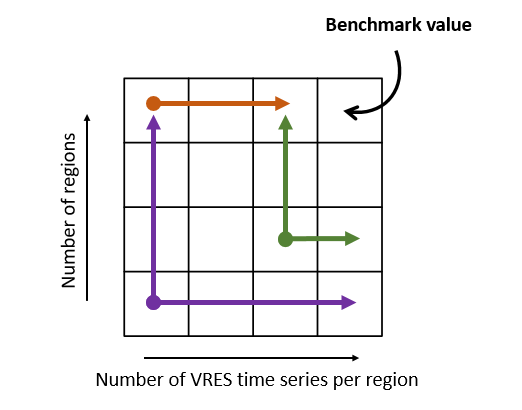}    
\caption{Pictorial depiction of how the optimal combinations of number of regions and number of VRES time series per region are determined from the TAC matrix. In order for a combination to be deemed optimal, the corresponding TAC and the TAC in the direction of increase of number of regions and number of VRES time series per region should not deviate more than $\pm5\%$ from the benchmark TAC.} 
\label{fig:tac_error_eg}
\end{center}
\end{figure}

With an aim to determine a set of optimal parameter combinations, a matrix containing TAC for each parameter combination is traversed. For each combination, the TACs in both the directions of increase are compared with the benchmark value as can be seen in Figure \ref{fig:tac_error_eg}. If the percentage error between the benchmark TAC and each of the candidate TACs is below an acceptable error threshold, then the parameter combination is deemed to be optimal. Here, an error threshold of $\pm5\%$ is assumed. In terms of absolute value, the error threshold is approximately $\pm10$ billion Euros/annum. 

Using this method, the optimal combinations obtained and the corresponding run times and TAC deviations are given in Table \ref{tab:optimalcombinations}. It is noteworthy that in the case of 33 regions, the TAC is under-estimated and it all other cases, the TAC is over-estimated. Among the optimal combinations, the combination of 33 regions and 38 VRES types is deemed to be the most optimal, as it is the lowest spatial resolution in the list.

\begin{table*}[h]
	\caption{Optimal combinations of parameters with corresponding TAC deviation and decrease in run time, compared to the benchmark}
	\begin{center}
	\begin{tabularx}{\textwidth}{c c c c c c}
		\toprule
		\thead{\textbf{Number of}\\\textbf{aggregated regions}} 
		& \thead{\textbf{Number of}\\\textbf{VRES time series}\\\textbf{per region}}  
		& \thead{\textbf{TAC deviation}\\\textbf{{[}1e9 Euro/annum{]}}}  
		& \thead{\textbf{TAC deviation}\\\textbf{{[}$\%${]}}} 
		& \thead{\textbf{Decrease in}\\\textbf{run time {[}min{]}}} 
		& \thead{\textbf{Decrease in}\\\textbf{run time {[}$\%${]}}}   \\
		\toprule
		
		\multirow{5}{*}{96}  & 58 & 0.55 & 0.26 & 131.62 & 28.07 \\
    		                & 48 & 1.29 & 0.62 & 83.49 & 17.81 \\
    		                & 38 & 2.8 & 1.33 & 143.21 & 30.55 \\
    		                & 28 & 4.18 & 1.99 & 165.43 & 35.29 \\
    		                & 18 & 6.39 & 3.04 & 189.61 & 40.44 \\
        \midrule
        \multirow{4}{*}{66}    & 68 & 1.56 & 0.74 & 250.69 & 53.47 \\ 
                        & 58 & 3.44 & 1.64 & 278.84 & 59.48 \\
                        & 48 & 4.69 & 2.23 & 305.31 & 65.12 \\
                        & 38 & 5.03 & 2.4 & 327.57 & 69.87 \\
        \midrule
        \multirow{4}{*}{50}    & 68 & 2.1 & 1.0 & 334.41 & 71.33 \\
                        & 58 & 2.64 & 1.26 & 343.18 & 73.2 \\
                        & 48 & 2.92 & 1.39 & 374.58 & 79.9 \\
                        & 38 & 4.1 & 1.95 & 403.67 & 86.1 \\
        
        \midrule
        \multirow{4}{*}{33}    & 68  & -5.23 & -2.49 & 421.99 & 90.01 \\
                         & 58  & -10.44 & -4.98 & 431.69 & 92.08 \\ 
                         & 48  & -4.75 & -2.26 &  436.7 & 93.15 \\ 
                          & 38  &  -9.29 & -4.42 & 435.79 & 92.95 \\ 
        
        \bottomrule
	\end{tabularx}
    \end{center}
	\label{tab:optimalcombinations}
\end{table*}

\section{Summary and Discussion} \label{sec:SummaryandDiscussion}

This study investigated the impact of spatial aggregation of regions and aggregation of VRES within each region, on optimal energy system design. For different combinations of number of regions and number of VRES within each region, the objective TAC resulting from model optimization was compared to the benchmark setting of 96 regions and 68 VRES types in each region. 

Further, since it is commonplace to consider administrative regions such as national borders while designing an optimal energy system, the study also considered 33 national regions. Also, it is commonplace to aggregate all system components including VRES within each model region to obtain just one representative type. This setting was also considered in the study. 

The TAC obtained in the case of 33 national regions and 33 aggregated regions were compared to the benchmark. It was observed that the TAC obtained in the case of 33 national regions deviated further from the benchmark, compared to the case of 33 aggregated regions. Since each model region is assumed to be a copper plate and the system components are aggregated in each region, it is important to pay attention to the region definitions. The defined regions should have similar component characteristics such that an aggregation of these components would not lead to loss of information, which in turn leads to inaccurate system design. 

To that end, considering very low number of aggregated regions (for e.g., 6 regions in Europe) would also result in inaccurate system design. In such cases, there is an over-simplification of the region network. In other words, the geographical gaps between generation sites and demand are largely ignored. This leads to under-estimation of the TAC. 

On the other hand, when 1 VRES type per region is considered, the temporal fluctuations present in the original set of VRES time series are not captured well. Further, the TAC is over-estimated because cost-competitive locations are not identified at such low spatial resolution of VRES. 

As the number of aggregated regions and VRES in each region was increased, the TAC deviation from benchmark was seen to reduce. At 33 aggregated regions and 38 representative VRES in each region, the system cost is under-estimated by $4.42\%$ and the run-time is reduced by $92.95\%$, compared to the benchmark. A further increase in the spatial resolution does not improve the results significantly, thereby deeming this setting to be the optimal.  

\section{Conclusion} \label{sec:Conclusion}
In this paper, a novel two-step aggregation scheme is introduced to reduce the underlying computational complexity of energy system optimization models. It involves two steps - (i) spatially aggregating homogeneous and spatially contiguous regions and (ii) aggregating variable renewable technologies such as wind turbines and photovoltaics in each newly-defined region. 

In order to aggregate the regions, we introduce a holistic approach that considers all the model parameters to find and group homogeneous regions, as opposed to considering just renewable potentials or transmission grid as is the common practice. Since each model region is assumed to be a copper plate, it is important to ensure that only regions that are spatially contiguous are grouped. Otherwise, the resulting regions are fragmented and are not fully connected. To that end, we employed the Hess model with additional contiguity constraints to ensure spatial contiguity in the resulting region groups. 

In order to reduce the number of wind turbines and photovoltaics in each region, we employed agglomerative hierarchical clustering to cluster each technology based on the similarity in their temporal profiles, to obtain a representative set. Such a method reduces the computational complexity while capturing the temporal fluctuations that are crucial for accurate energy system design. The common practice is to simply aggregate each technology in every region to obtain just one representative type. Our results show that keeping the number of wind turbine and photovoltaic types reasonably high in each region decreases the system costs because cost-competitive locations are identified. With 96 regions, if one technology type per region is considered the system cost is over-estimated by 30.11 billion Euro/annum, from the benchmark of 209.91 billion Euro/annum. 

\section{Outlook} \label{sec:Outlook}

\begin{figure*}[!htbp]
    \begin{center}
      \includegraphics[width=0.7\textwidth]{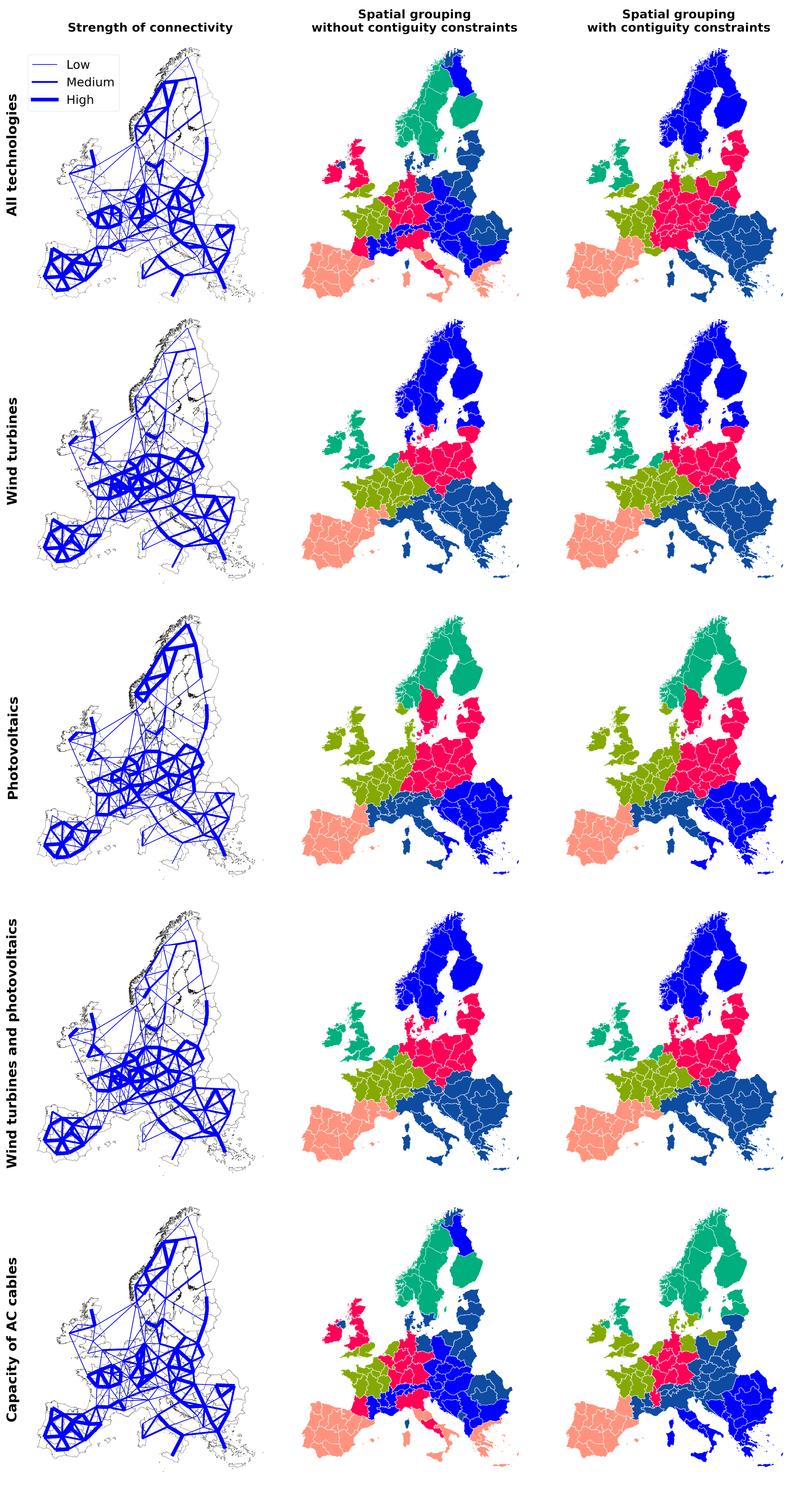}
      \caption{Region groups obtained when the 96 regions are spatially aggregated, to obtain 6 regions based on different technologies.}
        \label{fig:grouping_6reg_all}
    \end{center}
\end{figure*}

In this study, spatial aggregation of regions was performed based on all the components present in the energy system model. It would be interesting to benchmark the results against the aggregation based on specific components and their attributes. In Figure \ref{fig:grouping_6reg_all}, the region groups formed when the original regions are aggregated based on only all components, specific components, a combination of components, and a particular attribute of a component (capacity of AC cables), to obtain 6 regions is shown. Certainly, the region definitions in each case vary. However, it is interesting to see that when the regions are aggregated based on both wind turbines and photovoltaics, the regions look at lot similar to the ones obtained when aggregation is performed based on wind turbines. This shows that certain model components have bigger influence on the region definitions than others. Further, such an option to weight the components allows the modeller to choose the components that should influence the region definitions, depending on the research question. 

Our focus in this paper was to sufficiently represent variable renewable energy sources in each region. It would also be interesting to increase the spatial details of other components such as storage technologies and investigate the extent of its influence on energy system design. 

Clubbing the above-mentioned areas of future research, certain combinations of spatial and technology aggregations would also be worth investigating. For example, spatial aggregation based on electricity grid, and maintaining the spatial resolution of both source and storage technologies sufficiently high in each region. 

Finally, we have first performed spatial and technology aggregation and then a temporal aggregation in our experiments. It would be worth investigating the effect of a switch in this order, and to determine an optimal combination of spatial, technological, and temporal resolutions. 

\FloatBarrier 

\section{Code Availability} \label{sec:codeAvailibility}

The methods introduced in this paper are published in the Python package  \href{https://github.com/FZJ-IEK3-VSA/FINE}{FINE - Framework for Integrated Energy System Assessment} and can be easily applied and extended. 

\begin{ack}
The authors acknowledge the financial support of the Federal Ministry for Economic Affairs and Energy of Germany for the project METIS (project number 03ET4064).
\end{ack}

\bibliography{ifacconf}             

\end{document}